\documentclass{article}
\usepackage{amssymb,amsmath}
\usepackage{graphics}

\def\qed{\hfill {\hbox{${\vcenter{\vbox{               
   \hrule height 0.4pt\hbox{\vrule width 0.4pt height 6pt
   \kern5pt\vrule width 0.4pt}\hrule height 0.4pt}}}$}}}

\def\bar{\overline}

\newtheorem{theorem}{Theorem}
\newtheorem{definition}{Definition}
\newtheorem{lemma}[theorem]{Lemma}

\newtheorem{example}{Example}

\newenvironment{proof}[1][Proof]{\smallskip\noindent{\bf #1.}\quad}%
{\qed\par\medskip}

\author{
{\begin{tabular}{c} Daisy Lam \\
\small{\texttt{dlam002@student.ucr.edu}}\end{tabular}}
\and
{\begin{tabular}{c} Sam Nelson \\
\small{\texttt{knots@esotericka.org}}\end{tabular}}
\and
\small{\begin{tabular}{c}
Department of Mathematics, University of California, Riverside\\
900 University Avenue, Riverside, CA, 92521\end{tabular} }}

\date{}

\title{\Large \textbf{An isomorphism theorem for Alexander biquandles}}

\begin{document}
\maketitle

\begin{abstract}
We show that two Alexander biquandles $M$ and $M'$ are isomorphic iff 
there is an isomorphism of $\mathbb{Z}[s^{\pm 1},t^{\pm 1}]$-modules 
$h:(1-st)M\to (1-st)M'$ and a bijection $g:O_s(A)\to O_s(A')$ between the 
$s$-orbits of sets of coset representatives of $M/(1-st)M$ and $M'/(1-st)M'$ 
respectively satisfying certain compatibility conditions. 
\end{abstract}



\section{\large \textbf{Introduction}}

In \cite{N1}, it was shown that two finite Alexander quandles of the same
cardinality are isomorphic iff their $(1-t)$ submodules are isomorphic as
$\mathbb{Z}[t^{\pm 1}]$-modules. These finite Alexander
quandles are useful as a source of knot invariants defined by counting
homomorphisms in various ways -- setwise, weighted by cocyles in various 
quandle cohomology theories, etc. (see \cite{C2} for more).

Alexander quandles have been generalized to Alexander biquandles
\cite{K3}. Both quandles and biquandles are examples of algebraic structures 
with axioms derived from Reidemeister moves, the former with generators of
the algebra corresponding to arcs and the latter with generators corresponding
to semi-arcs in the knot diagram. The resulting non-associative algebraic 
structures are thus naturally suited for defining invariants of knots and 
links.

Biquandles have been studied in several recent papers such as \cite{FJK}, 
\cite{C1} and \cite{KM}. In particular, \cite{FJK} lists a number of known 
types of finite biquandles, including Alexander biquandles. In this paper we
give necessary and sufficient conditions for the existence of an isomorphism 
$f:M\to M'$ of Alexander biquandles which generalizes the main result from 
\cite{N1}.

\section{\large \textbf{Biquandles}}

We begin with the definition of a biquandle.

\begin{definition}
\textup{
A \textit{biquandle} is a set $B$ with four binary operations
$B\times B\to B$ denoted by
\[(a,b) \mapsto a^b, \ a^{\bar b}, \ a_b,\quad \mathrm{and} \quad a_{\bar b}\]
respectively, satisfying:}

\begin{list}{}{}
\item[\textup{1.}]{\textup{For every pair of elements $a,b\in B$, we have}
\[
\mathrm{(i)} \ a=a^{b{\bar{b_a}}}, \quad
\mathrm{(ii)} \ b=b_{a{\bar{a^b}}}, \quad
\mathrm{(iii)} \ a=a^{\bar{b}b_{\bar a}}, \quad
\mathrm{and} \quad
\mathrm{(iv)} \ b=b_{\bar{a}a^{\bar b}}.
\]}
\item[\textup{2.}]{\textup{ Given elements $a,b\in B$, there are unique
elements $x,y\in B$, possibly but not necessarily distinct, such that}
\[
\mathrm{(i)} \ x=a^{b_{\bar x}}, \quad
\mathrm{ (ii) } \ a=x^{\bar b},\quad
\mathrm{ (iii)} \   b=b_{{\bar x}a}, \]\[
\mathrm{ (iv)} \  y=a^{\bar{b_y}},\quad
\mathrm{ (v) } \  a=y^b,   \quad \mathrm{and} \quad
\mathrm{ (vi)} \  b=b_{y{\bar a}}.
\]}
\item[\textup{3.}]{\textup{ For every triple $a,b,c \in B$ we have:}
\[
\mathrm{ (i)} \ a^{bc}=a^{c_bb^c}, \quad
\mathrm{(ii)} \ c_{ba} =c_{a^bb_a}, \quad
\mathrm{(iii)} \ (b_a)^{c_{a^b}}=(b^c)_{a^{c_b}}, \]
\[
\mathrm{(iv)} \ a^{{\bar b}{\bar c}}
=a^{\bar{c_{\bar b}}\bar{b^{\bar c}}},
\quad
\mathrm{(v)} \ c_{{\bar b}{\bar a}} 
=c_{\bar{a^{\bar b}}{\bar {b_{\bar a}}}},
\quad \mathrm{and}  \quad
\mathrm{(vi)} \ (b_{\bar a})^{\bar{c_{{\bar {a^{\bar b}}}}}}
=(b^{\bar c})_{\bar{a^{\bar{c_{\bar b}}}}}.
\]}
\item[\textup{4.}]{\textup{Given an element $a\in B,$ there are 
unique elements
$x,y\in B$, possibly but not necessarily distinct, such that}
\[
\mathrm{ (i)} \   x=a_x, \quad
\mathrm{ (ii)} \  a=x^a, \quad
\mathrm{(iii)} \  y=a^{\bar y}, \quad \mathrm{and}  \quad
\mathrm{(iv) } \  a=y_{\bar a}.
\]
}
\end{list}
\end{definition}

A biquandle is a type of invertible \textit{switch}, i.e., a solution 
$S:X\times X\to X\times X$ to the (set-theoretic) \textit{Yang-Baxter 
equation}
\[(S\times I)(I \times S)(S\times I) = (I \times S)(S\times I)(I \times S)\]
where $X$ is a set and $I:X\to X$ is the identity. The components of such a 
solution $S$ satisfy axiom (3), and if $S$ is invertible the components and
the components of the inverse also satisfy axiom (1).
An invertible switch  $S(a,b)=(b_a,a^b)$ then defines a biquandle if its
component functions satisfy axioms (2) and (4). See \cite{FJK} for more.

The biquandle axioms are motivated by the Reidemeister moves in knot theory 
-- if we assign generators to each semi-arc in an oriented link diagram and 
consider the outbound semi-arcs at a crossing to be the results of the 
inbound semi-arcs operating on each other, with barred operations at negative 
crossings and unbarred operations at positive crossings, then the biquandle 
axioms are a set of minimal conditions required to make the resulting 
algebraic structure invariant under Reidemeister moves. 

In \cite{NV}, finite biquandles with cardinality $n$ are presented by 
$2n\times 2n$ block matrices composed of four $n\times n$ blocks which 
represent the operation tables of the four biquandle operations. 
Specifically, if $B=\{x_1,\dots, x_n\}$ then the matrix of $B$ has four 
blocks $M=\left[ \begin{array}{r|r} B^1 & B^2 \\ \hline B^3 & B^4 
\end{array}\right]$ such that
\[ B^l_{ij} =k \quad \mathrm{where} \quad  x_k=
\left\{\begin{array}{ll}
(x_i)^{(x_j)} & l=1 \\
(x_i)_{(x_j)} & l=2 \\
(x_i)^{\bar{(x_j)}} & l=3 \\
(x_i)_{\bar{(x_j)}} & l=4 \\
\end{array}
\right.\]

\begin{example}\textup{
The \textit{trivial} biquandle of order $n$ is the set $T=\{1,2,\dots, n\}$
with operations $i^j=i, \ i_j=i, \ i^{\bar j}=i$ and $i_{\bar j}=i$. It has
matrix}
\[
M=\left[\begin{array}{rrrr|rrrrr}
1 & 1 & \dots & 1 & 1 & 1 & \dots & 1 \\
2 & 2 & \dots & 2 & 2 & 2 & \dots & 2 \\
\vdots &  & & \vdots & \vdots & & & \vdots \\
n & n & \dots & n & n & n & \dots & n \\ \hline
1 & 1 & \dots & 1 & 1 & 1 & \dots & 1 \\
2 & 2 & \dots & 2 & 2 & 2 & \dots & 2 \\
\vdots & & & \vdots & \vdots & & & \vdots \\
n & n & \dots & n & n & n & \dots & n
\end{array}
\right].
\]
\end{example}

This matrix presentation was used in \cite{NV} to do a computer search in 
which all biquandles of order up to 4 were classified; matrix presentation 
of finite biquandles also makes symbolic computation with finite biquandles 
easy (see \cite{CN}). In \cite{FJK} and \cite{BF}, several examples of 
biquandle structures defined on groups and modules are given.
 
\begin{example} \textup{The following definition comes from \cite{BF}.
Let $M$ be a module over a ring $R$. Then 
$x^y=Cx+Dy$ and $x_y= Ay+Bx$ where $A,B\in R$ are invertible,
$C=A^{-1}B^{-1}A(I-A)$, and
$D= I-A^{-1}B^{-1}AB$ defines an invertible switch on $M\times M$ if the 
equation}
\[ [B,(A-I)(A,B)]=0\]
\textup{where $[X,Y]= XY-YX$ and $(X,Y)=X^{-1}Y^{-1}XY$ is satisfied. Thus,
module theory is a source of biquandles.}
\end{example}

\begin{example}
\textup{For a related example, let $M=\mathbb{Z}_2\oplus\mathbb{Z}_2$ 
considered as a $\mathbb{Z}_2$-module and set}
\[A=\left[\begin{array}{rr} 0 & 1 \\ 1 & 1 \end{array}\right], \quad 
B=\left[\begin{array}{rr} 1 & 1 \\ 0 & 1 \end{array}\right], \quad
C=\left[\begin{array}{rr} 1 & 0 \\ 1 & 1 \end{array}\right],\quad \mathrm{and}
\quad D=\left[\begin{array}{rr} 1 & 1 \\ 1 & 0 \end{array}\right].\] 
\textup{Then $M$ is a biquandle with}
\[x^y=x_{\bar{y}}=Cx+Dy+\left[\begin{array}{r} 1 \\ 1\end{array}\right] 
\quad \mathrm{and} \quad x_y=x^{\bar{y}}= 
Ay+Bx+\left[\begin{array}{r} 1 \\ 1\end{array}\right].\]
\textup{$M$ has biquandle matrix (where $x_1=
\left[\begin{array}{r} 1 \\ 0 \end{array}\right],\ 
x_2=\left[\begin{array}{r} 0 \\ 1 \end{array}\right],\ 
x_3=\left[\begin{array}{r} 1 \\ 1 \end{array}\right]$
and $x_4=\left[\begin{array}{r} 0 \\ 0 \end{array}\right]$)}
\[\left[\begin{array}{rrrr|rrrr}
3 & 1 & 2 & 4 & 4 & 1 & 3 & 2 \\
2 & 4 & 3 & 1 & 2 & 3 & 1 & 4 \\
1 & 3 & 4 & 2 & 3 & 2 & 4 & 1 \\
4 & 2 & 1 & 3 & 1 & 4 & 2 & 3 \\ \hline
4 & 1 & 3 & 2 & 3 & 1 & 2 & 4 \\
2 & 3 & 1 & 4 & 2 & 4 & 3 & 1 \\
3 & 2 & 4 & 1 & 1 & 3 & 4 & 2 \\
1 & 4 & 2 & 3 & 4 & 2 & 1 & 3 
\end{array}\right].\]

\textup{The counting invariant associated to this biquandle, 
$|\mathrm{Hom}(B(K),M)|$ where $B(K)$ is a knot biquandle, distinguishes
all of the Kishino knots from the unknot. See \cite{NV}.}
\end{example}

\section{\large \textbf{Alexander biquandles}}

In this section we give necessary and sufficient conditions for two Alexander
biquandles to be isomorphic. We begin with a definition from \cite{K3}.

\begin{definition} \label{d1}\textup{
Let $M$ be a module over the ring $\mathbb{Z}[s^{\pm 1},t^{\pm 1}]$ of 
Laurent polynomials in two variables. Then $M$ is a biquandle with operations}
\[ x^y=tx+(1-st) y, \quad x_y= sx, 
\quad x^{\bar y} = t^{-1}x + (1-s^{-1}t^{-1})y,
\quad \mathrm{and} \ x_{\bar y} = s^{-1} x.\]
\textup{Such a biquandle is called an \textit{Alexander biquandle}. }
\end{definition}

As expected, we have:

\begin{definition} \label{d2}\textup{
A 
\textit{homomorphism} of Alexander biquandles is a map $f:M\to M'$ satisfying}
\[f(x^y)=f(x)^{f(y)}, \ f(x_y)=f(x)_{f(y)}, \ f(x^{\bar y})=f(x)^{\bar{f(y)}},
 \ \mathrm{and} \ f(x_{\bar y})=f(x)_{\bar{f(y)}} \]
\textup{or equivalently}
\[f(tx+(1-st) y) =t f(x) + (1-st) f(y), \quad f(sx) =sf(x), \]\[
f(t^{-1}x+(1-s^{-1}t^{-1}) y) =t^{-1} f(x) + (1-s^{-1}t^{-1}) f(y), \quad 
\mathrm{and} \ f(s^{-1}x) =s^{-1}f(x).
\]
\end{definition}

\begin{example}\textup{If $s,t$ are invertible in $\mathbb{Z}_n$ then 
$\mathbb{Z}_n$ has the structure of an Alexander biquandle with the operations
above. For example, $\mathbb{Z}_3$ with $s=2$ and $t=1$ has $(1-st)=2$ and 
biquandle matrix}
\[
\left[\begin{array}{rrr|rrr}
3 & 2 & 1 & 3 & 2 & 1 \\ 
1 & 3 & 2 & 1 & 3 & 2 \\
2 & 1 & 3 & 2 & 1 & 3 \\ \hline
2 & 2 & 2 & 2 & 2 & 2 \\
1 & 1 & 1 & 1 & 1 & 1 \\
3 & 3 & 3 & 3 & 3 & 3
\end{array}\right]
\]
\textup{where we use $\mathbb{Z}_3=\{1,2,3\}$.}
\end{example}

If $f(0)=0$, then for $f:M\to M'$ to be a homomorphism of Alexander 
biquandles,
it suffices for $f$ to satisfy the first two equations in definition \ref{d2}:

\begin{lemma}\label{l3}
Let $f:M\to M'$ be a function which satisfies $f(0)=0\in M'$ and
\[f(tx+(1-st) y) =t f(x) + (1-st) f(y) \quad \mathrm{and} \quad 
f(sx) =sf(x). \] Then $f$ is a homomorphism of Alexander biquandles.
\end{lemma}

\begin{proof}
We must show that 
\[ f(t^{-1}x+(1-s^{-1}t^{-1}) y) =t^{-1} f(x) + (1-s^{-1}t^{-1}) f(y) \quad 
\mathrm{and} \ f(s^{-1}x) =s^{-1}f(x).\]
The second is easy: 
\[s^{-1}f(x)=s^{-1}f(s(s^{-1}x))=s^{-1}sf(s^{-1}x)=f(s^{-1}x).\]

The condition that $f(0)=0$ implies
\[f(tx)=f(tx+(1-st)0)=tf(x)+(1-st)f(0)=tf(x)\] and
\[f((1-st)y)=f(t0+(1-st)y)=tf(0)+(1-st)f(y)=(1-st)f(y).\] 
Then we also have
\[t^{-1}f(x)=t^{-1}f(t(t^{-1}x))=t^{-1}tf(t^{-1}x)=f(t^{-1}x).\]
Moreover, if $y=(1-st)z=tw$ then $f(-y)=-f(y)$ since
\begin{eqnarray*} f(-y)+f(y) & = & f(t(-w))+f((1-st)z) \\
& = & tf(-w)+(1-st)f(z) \\
& = & f(t(-w)+(1-st)z) \\
& = & f(-y+y)=f(0)=0.\end{eqnarray*}
But then
\begin{eqnarray*}
f(t^{-1}x+(1-s^{-1}t^{-1})y) & = & f(t(t^{-2}x) +(1-st)(-s^{-1}t^{-1}y)) \\
& = & tf(t^{-2}x)+(1-st)f(-s^{-1}t^{-1}y) \\
& = & t^{-1}f(x) +f(-s^{-1}t^{-1}(1-st)y) \\
& = & t^{-1}f(x) -s^{-1}t^{-1}f((1-st)y) \\
& = & t^{-1}f(x) +(1-st)(-s^{-1}t^{-1})f(y) \\
& = & t^{-1}f(x) +(1-s^{-1}t^{-1})f(y)
\end{eqnarray*} 
as required.
\end{proof}

\begin{lemma} \label{l1}
Let $f:M\to M'$ be a homomorphism of biquandles. Then $f(0)\in Ker(\phi)$
where $\phi:M'\to M'$ is given by $\phi(x)=(1-s)x$.
\end{lemma}

\begin{proof}
Since $f$ is a homomorphism of biquandles, we have $f(x_y)=f(x)_{f(y)}$ for
all $x,y\in M$. In particular, $f(0)=f(s0)=f(0_y)=f(0)_{f(y)}=sf(0)$
and we have $(1-s)f(0)=0\in M'.$
\end{proof}

\begin{lemma} \label{l2}
Let $g_z:M'\to M'$ be defined by $g_z(x)=x+z$. Then $g_z$ is an isomorphism
of Alexander biquandles if $z\in \mathrm{Ker}(\phi)$ where $\phi:M'\to M'$ 
is given by $\phi(x)=(1-s)x$.
\end{lemma}

\begin{proof}
For any $z\in M'$, $g_z(x)=x+z$ is bijective. Thus,
we must show that $(1-s)z=0$ implies that $g_z$ is a homomorphism of
biquandles. That is, we must compare
\begin{list}{}{}
\item{(1)}{\quad $g_z(ta+(1-st)b)$ with $tg_z(a)+(1-st)g_z(b)$,}
\item{(2)}{\quad $g_z(sa)$ with $sg_z(a),$}
\item{(3)}{\quad $g_z(t^{-1}a+(1-s^{-1}t^{-1})b)$ with
$t^{-1}g_z(a)+(1-s^{-1}t^{-1})g_z(b)$ and}
\item{(4)}{\quad $g_z(s^{-1}a)$ with $s^{-1}g_z(a)$}
\end{list}
where $a,b\in M'$.

For (1) we see that
\[
g_z(ta+(1-st)b)  =  ta+(1-st)b +z \]
and \[ tg_z(a)+(1-st)g_z(b)  =  ta + tz +(1-st) b +(1-st)z
\]
so subtracting yields
\begin{eqnarray*}
-ta-(1-st)b-z+ta+tz+(1-st)b+(1-st)z & = & -z+tz+(1-st)z \\
 & = & -z+tz+z-stz \\
 & = & t(1-s)z =0.
\end{eqnarray*}

For (2) we see that
\[ g_z(sa)  =  sa+z \quad \mathrm{and} \quad sg_z(a)=s(a+z) \]
so subtracting yields
\begin{eqnarray*}
sa+z -s(a+z) & = & sa+z -sa -sz \\
& = & z-sz =(1-s) z =0.
\end{eqnarray*}

Finally, by lemma \ref{l3}, we are done.

\end{proof}

Not every biquandle isomorphism $f:M\to M'$ sends $0\in M$ to $0\in M'$,
but in light of lemmas \ref{l1} and \ref{l2}, we may replace any isomorphism
$f:M\to M'$ which does not with $f'=g_{(-f(0))}\circ f$, and then $f'(0)=0$.

Let us denote the orbit of a subset $X\subseteq M$ under multiplication by 
$s$ by
\[O_s(X)=\{ s^ix \ | \ i\in \mathbb{Z}, x\in X \}.\]

We can now state our main result.

\begin{theorem} \label{tmain}
Two Alexander biquandles $M$ and $M'$ are isomorphic as biquandles iff they
satisfy 
\begin{list}{}{}
\item[(i)]{There is an isomorphism of
 $\mathbb{Z}[s^{\pm 1},t^{\pm 1}]$-modules $h:(1-st)M\to(1-st)M'$ and}
\item[(ii)]{For every set of coset representatives $A$ of $M/(1-st)M$ 
in which the class of $(1-st)M$ is represented by $0\in M$, there
is a corresponding set of coset representatives $A'$ of $M'/(1-st)M'$ and a
bijection $g:O_s(A) \to O_s(A')$ such that $g(A)=A'$ and
\[(1-st)g(\alpha)=h((1-st)\alpha) \quad \mathrm{and} \quad 
g(s\alpha+\omega) = sg(\alpha)+h(\omega)\] for 
every $\alpha\in A$ and $\omega\in (1-st)M$.}
\end{list}
\end{theorem}

\begin{proof}
\noindent($\Rightarrow$) Suppose $f:M\to M'$ is an isomorphism of biquandles, 
and without loss of generality suppose $f(0)=0.$ Then $f$ commutes with 
multiplication by powers of $s$, $t$ and $(1-st)$ and satisfies
\[f(x+y)=f(x)+f(y) \quad \quad \quad \forall x,y\in (1-st)M\] 
since $x,y\in (1-st)M$ implies $x=tx', y=(1-st)y'$ and then
\begin{eqnarray*}
f(x+y)& = & f(tx'+(1-st)y')\\
& = & tf(x')+(1-st)f(y')\\
& = & f(tx')+f((1-st)y') \\
& = & f(x)+f(y).
\end{eqnarray*}
Hence, $h=f|_{(1-st)M}$ is an isomorphism of 
$\mathbb{Z}[s^{\pm 1},t^{\pm 1}]$-modules.

Now, let $A=\{\alpha_i \ | \ i\in I\}$ for some indexing set $I$ be a set 
of coset 
representatives for $M/(1-st)M$ and define $g=f|_{O_s(A)}$. The image of
$g$ is $f(O_s(A))=O_s(f(A))$ since $f$ commutes with $s$, and thus 
$g$ is bijective.
Then every element of $M$ has the form $x=\alpha_i+(1-st)\omega$ for some
$\alpha_i\in A$ and $(1-st)\omega\in (1-st)M$, and we have
\begin{eqnarray*}
f(x)& = & f(tt^{-1}\alpha_i+(1-st)\omega)\\
& = & tf(t^{-1}\alpha_i)+(1-st)f(\omega) \\
& = & f(tt^{-1}\alpha_i)+f((1-st)\omega) \\
& = & g(\alpha_i)+h((1-st)\omega).
\end{eqnarray*}
In particular,
\[f(\alpha_i+(1-st)M)=g(\alpha_i)+h((1-st)M)=g(\alpha_i)+(1-st)M';\]
that $A'=\{g(\alpha_i)\ | \ i\in I\}$ is a set of coset 
representatives of $M'/(1-st)M'$ then follows from the bijectivity of
$f$ and the 
fact that if $g(\alpha_i)=g(\alpha_j)+(1-st)m'$, then
\begin{eqnarray*}
\alpha_i & = & f^{-1}f(\alpha_i) \\
& = & f^{-1}(g(\alpha_j)+(1-st)m') \\
& = & \alpha_j+f^{-1}(1-st)m'
\end{eqnarray*} and hence $\alpha_i=\alpha_j$.

For any $\alpha\in A$ we have
\begin{eqnarray*}
(1-st)g(\alpha) & = & 0+(1-st)f(\alpha) \\
& = & tf(0)+(1-st)f(\alpha) \\
& = & f(t0+(1-st)\alpha) \\
& = & f((1-st)\alpha)\\
& = & h((1-st)\alpha).
\end{eqnarray*} 

Moreover, if $s\alpha+\omega\in O_s(A)$ where $\alpha\in A$ and 
$\omega=(1-st)\omega'\in (1-st)M$ then
\begin{eqnarray*}
g(s\alpha+\omega) & = & f(s\alpha+\omega) \\
& = & f(tst^{-1}\alpha+(1-st)\omega') \\
& = & tf(st^{-1}\alpha)+(1-st)f(\omega')\\
& = & tf(st^{-1}\alpha)+f((1-st)\omega')\\
& = & sf(\alpha) +f(\omega)\\
& = & sg(\alpha)+h(\omega)
\end{eqnarray*}
as required.

\smallskip

\noindent
($\Leftarrow$) Suppose that $h:(1-st)M\to (1-st)M'$ is an isomorphism 
of $\mathbb{Z}[s^{\pm 1},t^{\pm 1}]$-modules, $A$ and $A'$ 
respectively are sets of coset representatives for $M/(1-st)M$ and
$M'/(1-st)M'$ with
$0\in M$ as the representative in $A$ of the coset $0+(1-st)M \subset M$,
and that $g:O_s(A)\to O_s(A')$ is a 
bijection satisfying $g(A)=A'$,
\[(1-st)g(\alpha)=h((1-st)\alpha)\quad \mathrm{and}\quad
g(s\alpha+\omega)=sg(\alpha)+h(\omega)\]
for all $\alpha\in A$ and $s\alpha+\omega\in O_s(A)$ with $\omega\in (1-st)M$.
In particular,
\[(1-st)g(0)=h((1-st)0)=h(0)=0\]
implies that $g(0)=0\in A'$.

Now define $f:M\to M'$ by setting
\[f(s^i\alpha+\omega)=g(s^i\alpha)+h(\omega) 
= s^ig(\alpha)+h(\omega), \quad \alpha\in A, \ \omega\in 
(1-st)M.\] To see that $f$ is well-defined, note that every element of
$M$ can be written as $x=\alpha+\omega$ in a unique way with $\alpha\in A$, 
$\omega\in (1-st)M$. So, if $x=\alpha+\omega=s^j\alpha'+\omega'+\omega''$ 
with $s^j\alpha'+\omega'\in O_s(A)$, then 
$s^j\alpha'+\omega'=\alpha+\omega-\omega'' \in O_s(A)$ and we have
\[g(s^j\alpha'+\omega')=g(\alpha+\omega-\omega'')
=g(\alpha)+h(\omega-\omega')\]
so that
\begin{eqnarray*}
g(s^j\alpha'+\omega')+h(\omega'') 
& = & g(\alpha)+h(\omega-\omega'') + h(\omega'') \\
& = & g(\alpha)+h(\omega)-h(\omega'') + h(\omega'') \\
& = & g(\alpha)+h(\omega).
\end{eqnarray*}
Define $k:A\to A'$ by $k=g|_{A}$. Then $k$ is bijective, and 
$f(\alpha+\omega)=k(\alpha)+h(\omega)$ for every $\alpha\in A,\omega\in 
(1-st)M$. Then $f$ is bijective, since $f$ is setwise the cartesian 
product of the bijective maps $k$ and $h$.

Now if $x=\alpha+\omega$ we have
\[f(sx)=f(s\alpha +s\omega) =g(s\alpha)+h(s\omega) =sg(\alpha)+sh(\omega)
=s(g(\alpha)+h(\omega))=sf(x).\] 
It follows that $s^ig(x)=g(s^ix)$ for every $i\in \mathbb{Z}.$

Note that $t\alpha = s^{-1}\alpha-s^{-1}(1-st)\alpha$ 
for any $\alpha\in M$. Then if $x=\alpha_1+\omega_1$ and $y=\alpha_2+\omega_2$
with $\alpha_i\in A$ and $\omega_i\in (1-st)M$, we have
\begin{eqnarray*}
f(tx+(1-st)y) & = &
f(t(\alpha_1+\omega_1) +(1-st)(\alpha_2+\omega_2)) \\ 
& = & f(t\alpha_1+t\omega_1 +(1-st)\alpha_2+(1-st)\omega_2) \\
& = & f(s^{-1}\alpha_1-s^{-1}(1-st)\alpha_1+t\omega_1 
+(1-st)\alpha_2+(1-st)\omega_2) \\
& = & g(s^{-1}\alpha_1) +h(-s^{-1}(1-st)\alpha_1+t\omega_1 
+(1-st)\alpha_2+(1-st)\omega_2) \\
& = & s^{-1}g(\alpha_1) -s^{-1}h((1-st)\alpha_1) + th(\omega_1) +
h((1-st)\alpha_2)+(1-st)h(\omega_2) \\
& = & s^{-1}g(\alpha_1) -s^{-1}(1-st)g(\alpha_1) + th(\omega_1) 
h((1-st)\alpha_2)+(1-st)h(\omega_2) \\
& = & tg(\alpha_1)+ th(\omega_1) + (1-st)g(\alpha_2)+(1-st)h(\omega_2) \\
& = & tf(x)+(1-st)f(y) \\
\end{eqnarray*}
and $f$ is an isomorphism of biquandles.
\end{proof}

If $s=1$ then $O_s(A)=A$ and $s\alpha+\omega\in O_s(A)$ implies $\omega=0$,
so the condition that $g(s\alpha+\omega)=sg(\omega)+h(\omega)$ is automatic.
It is then possible to show (see \cite{N1}) that if $|M|=|M'|<\infty$ and 
condition (i) is satisfied, then for every choice of coset representatives 
$A$ of $M/(1-st)M$ there exists a corresponding set of coset representatives 
$A'$ of $M'/(1-st)M'$ so that $(1-st)g(\alpha)=h((1-st)\alpha)$.
If $s\ne  1$, that is, if $M$ and $M'$ are Alexander biquandles which are not 
Alexander quandles, then condition (i) and $|M|=|M'|$ are not sufficient for
$M$ to be isomorphic to $M'$ as biquandles, as the next example demonstrates.

\begin{example}\textup{
Let $M=\mathbb{Z}_8$ with $s=3$ and $t=5$, and let $M'=\mathbb{Z}_8$
with $s=5$ and $t=3$. Then $(1-st)=6$ and $(1-st)M=(1-st)M'=\{0,2,4,6\}$.
Moreover, $h:(1-st)M\to(1-st)M'$ defined by $h(x)=-x$ satisfies 
$h(3x)=5h(x)$ and $h(5x)=3h(x)$ since $5(2)\equiv -3(2) \ \mathrm{mod} \ 8$.}

\textup{
Now, let $A=\{0,1\}$ so that $O_s(A)=\{0,1,3\}$. Then in order to satisfy}
\[6g(1)=h((6)1)=h(6)=2\] 
\textup{we must have either $g(1)=3$ or $g(1)=7$. If $g(1)=3$ then }
\[g(sx)=sg(x) \quad \Rightarrow \quad
g(3)=g(3(1))=5g(1)=5(3)=7\] \textup{but then} 
\[g(\alpha+\omega)=g(\alpha)+h(\omega) \quad \Rightarrow \quad
g(3)=g(1+2)=g(1)+h(2)=3+6=1\ne 7.\]
\textup{Similarly, $g(1)=7$ implies
$g(3(1))=5g(1)=5(7)=3$ while $g(1+2)=g(1)+h(2)=7+6=5\ne 3$. Thus, for our 
choice of coset representatives $A$ there is no bijection 
$g:O_s(A)\to O_s(A')$
satisfying $(1-st)g(\alpha)=h((1-st)\alpha)$ and 
$g(s\alpha+\omega)=sg(\alpha)+h(\omega)$ for all $\alpha, \alpha+\omega\in A$
and $\omega\in (1-st)M$, and hence $M$ and $M'$ are non-isomorphic 
Alexander biquandles. Our \textit{Maple} computations confirm this result.}
\end{example}

\section{\large \textbf{Future Research}}

There remains much to be done in the study of biquandles and Alexander 
biquandles in particular. Computation of the Yang-Baxter homology groups for
Alexander biquandles might shed additional light on the homology of quandles. 
It is conjectured that the knot biquandle is a complete invariant of virtual 
knots up to vertical mirror image (see \cite{FJK} and \cite{K4}); if this is 
true, then it 
should be possible to derive nearly all other invariants of knots from the 
biquandle of a knot, potentially illuminating the relationships between 
the invariants in the process.

\section{\large \textbf{Acknowledgements}}

We would like to thank the referee, whose helpful observations improved
and strengthened this paper.

\end{document}